\documentclass[11pt]{article}

\usepackage[utf8]{inputenc}
\usepackage{amssymb}
\usepackage{amsmath}

\input{amssym.def}
\input{amssym}

\newtheorem{theorem}{Theorem}

\newtheorem{definition}[theorem]{Definition}

\newtheorem{lemma}[theorem]{Lemma}
\newtheorem{proposition}[theorem]{Proposition}
\newtheorem{remark}[theorem]{Remark}

\def\de{{\delta}}

\def\pa{{\partial}}

\def\tL{\tilde{L}}
\def\ov{\overline}
\def\Y{\mathbf{Y}}

\def\qq{q^{-1}}

\def\lhq{{\cal L}(\h, q)}

\def\h{\hbar}

\def\C{{\Bbb C}}

\def\vv{V^{\otimes 2}}

\def\lhq{\ifmmode {\cal L}(q,\hbar)\else ${\cal L}(q,\hbar)$\fi}

\def\lqh{\ifmmode {\cal L}(q,\hbar)\else ${\cal L}(q,\hbar)$\fi}

\def\Tr{{\rm Tr}}

\def\al{{\alpha}}

\def\be{\begin{equation}}
\def\ee{\end{equation}}

\begin{document}

\makeatletter
\renewcommand{\theequation}{{\thesection}.{\arabic{equation}}}
\@addtoreset{equation}{section} \makeatother

\title{Bethe subalgebras in Braided  Yangians and Gaudin-type models}

\author{\rule{0pt}{7mm} Dimitri
Gurevich\thanks{gurevich@ihes.fr}\\
{\small\it Universit\'e de Valenciennes, EA 4015-LAMAV}\\
{\small\it F-59313 Valenciennes, France}\\
{\small\it and}\\
{\small\it Interdisciplinary Scientific Center Poncelet (ISCP, UMI 2615 du CNRS)}\\
{\small\it Moscow 119002, Russian Federation}\\
\rule{0pt}{7mm} Pavel Saponov\thanks{Pavel.Saponov@ihep.ru}\\
{\small\it
National Research University Higher School of Economics,}\\
{\small\it 20 Myasnitskaya Ulitsa, Moscow 101000, Russian Federation}\\
{\small \it and}\\
{\small \it
Institute for High Energy Physics, NRC "Kurchatov Institute"}\\
{\small \it Protvino 142281, Russian Federation}\\
\rule{0pt}{7mm} Alexei Slinkin\\
{\small\it
National Research University Higher School of Economics,}\\
{\small\it 20 Myasnitskaya Ulitsa, Moscow 101000, Russian Federation}\\
}

\maketitle

\begin{abstract} In \cite{GS1} the notion of braided Yangians of Reflection Equation type was introduced. Each of these algebras is associated with an involutive
or Hecke symmetry $R$. Besides, the quantum analogs of certain symmetric polynomials (elementary symmetric ones, power sums) were suggested.
In the present paper we show that these quantum symmetric polynomials  commute with each other and consequently  generate a commutative Bethe subalgebra.
As an application,  we get some Gaudin-type models and  the corresponding Bethe subalgebras.

\end{abstract}

{\bf AMS Mathematics Subject Classification, 2010:} 81R50, 82B23


\section{Introduction}

Let $V$ be a finite dimensional vector space over the ground field $\C$, $\dim_{\,\Bbb C}V = N$. Consider a linear operator $R\in{\rm End}(\vv)$ which is a solution
of the so-called braid relation
\be
R_1 R_2 R_1=R_2R_1 R_2. \label{braid} \ee
Both sides of this braid relation belong to ${\rm End}(V^{\otimes 3})$. Also, we use the standard shorthand  notation $R_k=R_{k\, k+1}$
for the operator $R$ acting in the positions $k$ and $k+1$ of the tensor product $V^p$, $p\geq k+1$.

An operator subject to (\ref{braid}) will be called a {\it braiding}.

In addition, we  assume all braidings $R$, we are dealing with, to obey either  the Hecke condition
$$
(R - q\, I)(R+ \qq\, I)=0,\quad q\in \C\setminus \{0,\pm 1\}
$$
or the condition
$$
R^2 = I.
$$
Hereafter, $I$ is the identity operator or matrix. Respectively, we will call such an operator $R$ a {\it  Hecke symmetry} or an {\it involutive symmetry}.

Given a Hecke or involutive symmetry $R$, we can construct some current (i.e. depending on parameters) operators
$$
R(u,v)=R-\frac{(q-\qq) \, u\, I}{u-v},
$$
provided $R$ is a Hecke symmetry,  and
$$
R(u,v)=R-\frac{I}{u-v}
$$
provided $R$ is involutive. These current operators satisfy the Quantum Yang-Baxter equation with parameters
\be
R_{1}(u,v)R_{2}(u,w)R_{1}(v,w)=R_{2}(v,w)R_{1}(u,w)R_{2}(u,v).
\label{YB}
\ee
The correspondence $R\mapsto R(u,v)$ is called the {\em Baxterization procedure} and the resulting operator $R(u,v)$ is called a trigonometric (resp., rational) $R$-matrix.

Since, in fact, the current operator $R(u,v)$ depends  on the ratio of parameters $x=u/v$ in the trigonometric case or on their difference $x=u-v$ in the rational case,
we also use the notation $R(x)$ for both cases. So, we respectively have
\be
R(x) = R - \frac{(q-q^{-1})x}{x-1}\,I,\qquad  R(x) = R-\frac{1}{x}\,I.
\label{R-trig}
\ee
We hope this does not lead to a misunderstanding.

In \cite{GS1, GS3} the so-called {\em braided Yangian} was introduced. It is  generated by an $N\times N$ current matrix $L(u)$, which is
subject to the system
\be
R_1(u,v)\,L_{\ov 1}(u)\, L_{\ov 2}(v)-L_{\ov 1}(v)\,L_{\ov 2}(u)\, R_1(u,v)=0.
\label{sys1}
\ee
Hereafter, we use the following notation
\be
L_{\ov 1}(u)=L_{1}(u),\quad L_{\ov{k+1}}(u)=R_{k}L_{\ov k}(u) R_{k}^{-1},\quad  k\ge 1.
\label{transf}
\ee

The generating matrix $L(u)$ is assumed to be a formal series in $u$
\be
L(u)=I+\sum_{k=1}^\infty L[k] \, u^{-k},
\label{gen-matr}
\ee
where the matrices $L[k]$, $k\ge 1$, are called the Laurent coefficients.
Consequently, the condition (\ref{sys1}) leads to an infinite system of relations on the Laurent coefficients $L[k]$. We use the notation $\Y(R)$ for the braided Yangian
defined by the system (\ref{sys1}). Below, we mainly deal with the whole matrix $L(u)$ without using its Laurent coefficients $L[k]$.

Note, that if we set $R=P$ (hereafter, $P$ stands for the usual flip) the system (\ref{sys1}) gives rise to the  Drinfeld's Yangians $\Y(gl(N))$.
By the initial definition, $\Y(gl(N))$ is an algebra, generated by entries of the corresponding  Laurent coefficients $L[k]$, whereas  entries of the matrix $L(u)$ are treated
to be elements of $\Y(gl(N))[[u^{-1}]]$.

\begin{remark}\rm
A Yangian-like algebra can be defined in a more general way which is based on a compatible pair of braidings $(R,F)$.
Recall that according to \cite{IOP} an ordered pair of braidings $(R,F)$ is called {\em compatible} if they are subject to the relations
\be
R_{1} F_{2} F_{1} = F_{2} F_{1}R_{2},\qquad R_{2} F_{1} F_{2} = F_{1} F_{2} R_{1}.
\label{coor}
\ee
Given such a pair, we define the generalized braided Yangian by the relation (\ref{sys1}) where the copies $L_{\overline k}(u)$ are defined with the use of
the braiding $F$:
$$
L_{\ov 1}(u)=L_{1}(u),\quad L_{\ov{k+1}}(u)=F_{k} L_{\ov k}(u) F_{k}^{-1},\quad  k\ge 1.
$$
For any compatible braidings $R$ and $F$ the use of (\ref{coor}) allows us to shift the relations (\ref{sys1}) to higher positions:
\be
R_k(u,v)L_{\overline k}(u)L_{\overline {k+1}}(v) = L_{\overline k}(v)L_{\overline {k+1}}(u) R_k(u,v), \quad \forall k\ge 1.
\label{sys1-sh}
\ee
Namely this property is the main ``raison d'\^{etre}" for the notion of compatible braidings.

It is easy to see that if $F=P$ or $F=R$, we get a compatible pair $(R,F)$ for any braiding $R$. In order to distinguish these specific cases we call
 the algebra $\Y(R)$ the {\em braided Yangian of Reflection Equation (RE) type}, while for $F=P$ we call the corresponding algebra the {\em Yangian of RTT type}.
If $R$ is the Hecke symmetry coming from the quantum group $U_q(sl(N))$, the corresponding Yangian of RTT type is somewhere called the {\it q-Yangian} (see \cite{M}).
\end{remark}

Note, that analogs of certain symmetric polynomials and the corresponding Bethe subalgebras as well as some analogs of the Newton and Cayley-Hamilton identities can be
constructed in all braided Yangians (not only in these of RE type).  However, the braided Yangians of RE type have a lot of specific properties, which distinguish them from
other braided Yangians. Thus, as was shown in \cite{GS1}, the evaluation morphisms in the braided Yangians of RE type are more similar to these in the  Drinfeld's Yangians
$\Y(gl(N))$.  Moreover,  the quantum determinant (the higher order elementary symmetric polynomial) is always central in any braided Yangian of RE type, while for other
types of braided Yangians it is not true in general (see \cite{GS1}).

In the present paper we show that if $R$ is an even symmetry, then the elementary symmetric polynomials in the corresponding  braided Yangian of RE type commute with
each other and consequently they generate a commutative subalgebra called {\it the Bethe one}. Also, we define some  Gaudin-type models and by applying the Talalaev's
method we construct their Bethe subalgebras, provided $R$ is an involutive symmetry. We call these models {\em braided Gaudin ones}.

The paper is organized as follows. The Bethe subalgebras in the braided Yangians are constructed in the next section. Their commutativity is established  in section 3. In section
4 we introduce  braided  versions of the rational Gaudin model and exhibit their Bethe subalgebras.

{\bf Acknowledgement}  The work of P.S. has been funded by the Russian Academic Excellence Project '5-100' and was also partially supported
by the RFBR grant 16-01-00562.

\section{Bethe subalgebras via quantum symmetric polynomials}

We begin this section with the following observation.
The current $R$-matrices, we are dealing with, are non-singular except for some isolated values of parameters. Namely, in the rational and trigonometric cases we respectively have
$$
R^{-1}(x) = \frac{x^2}{x^2-1}\,R(-x),\quad x\not=\pm 1,
$$

$$
R^{-1}(x) = \frac{(x-1)^2}{(x-1)^2-\lambda^2 x}\,R(-x^{-1}),\quad x\not=q^{\pm 2},\quad \lambda \equiv q-q^{-1}.
$$

Let us define a family of skew-symmetrizers $A^{(k)}_{1\dots k}(R)\in {\rm End}(V^{\otimes n})$, $n\ge k$, associated with a Hecke symmetriy $R$.
They are defined by the following recurrent relations:
\be
A^{(1)}=I,\qquad
A^{(k+1)}_{1\dots k+1} = \frac{k_q}{(k+1)_q}\,A^{(k)}_{1\dots k}\Big(\frac{q^k}{k_q} \,I - R_k\Big)A^{(k)}_{1\dots k},\quad k\ge 1,
\label{a-sym}
\ee
where $I$ is the identity matrix  of an appropriate size and  $k_q = (q^k - q^{-k})/(q-q^{-1}) $ is the standard notation for $q$-numbers.
The parameter $q$ is assumed to be generic: $q^k\not=1,\,\,\,\,\forall k\in {\Bbb Z}_+$, so all  $q$-integers are nonzero.  The low indeces indicate the positions where these operators are applied.

For an involutive symmetry $R$ the skew-symmetrizers are defined by the same formula (\ref{a-sym}) where we should set $q=1$ and use $k$ instead of $k_q$.

Hereafter, we confine ourselves to even skew-invertible  symmetries. The  bi-rank of these symmetries (Hecke or involutive) are always assumed to be
$(m|0)$\footnote{Observe that in general $m\not= N$. However, if $R$ equals $P$ or is its deformation, then $m=N$.}. A definition of even symmetries,  bi-rank, and a review of
properties of the corresponding symmetries can be found, for example, in \cite{GPS, GS1}. Here, we exhibit some of them:
\begin{itemize}
\item
For any  braiding under consideration there exists an $N\times N$ matrix $C$ with the following properties:
\be
{\rm Tr}_{(2)}R_{12}C_2 = I_1,\qquad R_kC_kC_{k+1} = C_kC_{k+1}R_k\quad \forall\,k\ge 1,\qquad {\rm Tr}\,C = \frac{m_q}{q^m}.
\label{C}
\ee
\item
With the matrix $C$ we define the {\it $R$-trace} of a matrix $X$ (of an appropriate size) with entries from an algebra $\cal A$:
\be
{\rm Tr}_{R(12\dots k)}X = {\rm Tr}_{(12\dots k)}(C_1C_2\dots C_k\,X), \quad X\in ({\rm Mat}_N({\cal A}))^{\otimes k}.
\label{R-tr}
\ee
The second relation in (\ref{C}) allows us to prove an important {\it cyclic property} of the $R$-trace:
\be
{\rm Tr}_{R(12\dots k)}(f(R_1,R_2\dots R_{k-1}) X) = {\rm Tr}_{R(12\dots k)}(X f(R_1,R_2\dots R_{k-1}))
\label{cyclic-p}
\ee
for any polynomial $f$ in $R$-matrices and for a matrix $X$ as above.
\item
The following property of the $R$-trace is of the great importance in what follows:
\be
{\rm Tr}_{R(k+1)}(R_k^{\pm 1}X_{1\dots k}R_{k}^{\mp 1}) = I_{1\dots k}{\rm Tr}_{R(k)}(X_{1\dots k}),
\label{tr-property}
\ee
where $X_{1\dots k}$ is an operator (a matrix) in the space $V^{\otimes k}$ and $I_{1\dots k}$ is the identity operator (matrix) in the same space.
\item
The skew-symmetrizers (\ref{a-sym}) are idempotents
\be
A^{(k)}A^{(k)} = A^{(k)}\quad \forall\, k\ge 1,
\label{a-sym-pro}
\ee
and they vanish for $k\geq m+1$:
$$
A^{(m+1)}\equiv 0, \qquad {\rm rank }\,A^{(m)} = 1.
$$

\end{itemize}

Now,  we introduce  quantum analogs of the elementary symmetric polynomials defined via  generating matrix $L(u)$ of a braided
Yangian of RE type $\Y(R)$ in the trigonometrical case.
\begin{definition}
\rm
The elements
\be
e_k (u) = {\rm Tr}_{R(12 \ldots k)} \,A^{(k)}_{12 \ldots k} {L}_{\overline{1}} (u) {L}_{\overline{2}} (q^{-2} u) \ldots {L}_{\overline{k}} (q^{-2(k-1)} u), \quad k\ge 1,
\label{elem-sym}
\ee
are called (quantum) elementary symmetric polynomials in the Yangian of RE type $\Y(R)$.

The subalgebra generated in $\Y(R)$  by the elementary symmetric polynomials is called the {\em Bethe subalgebra}.
\end{definition}

It should be emphasized that this definition (as well as the definition below of the power sums) is valid without assuming $R$ to be even. 
However, if $R$ is even of bi-rank  $(m|0)$, then there exists the highest elementary symmetric polynomial $e_m(u)$
which is called the {\em quantum determinant}. Observe that it is  always central in the corresponding braided Yangian $\Y(R)$.
 
In the next section we show that  the elementary symmetric polynomials commute with each other and consequently, the Bethe subalgebra is commutative.
To this end we need the following lemma.
\begin{lemma}
The following matrix identity holds true:\rm
\be
{\rm Tr}_{R(k+1\dots k+p)}A^{(p)}_{k+1\dots k+p}L_{\overline{k+1}}(u)L_{\overline{k+2}}(q^{-2}u)\dots L_{\overline{k+p}}(q^{-2(p-1)}u) = I_{12\dots k}\,e_p(u),
\label{elem-sym-sh}
\ee\it
where  $e_p(u)$ is an elementary symmetric polynomial defined in (\ref{elem-sym}).
\end{lemma}

\medskip

\noindent
{\bf Proof.}
First of all, note that due to the definition of $L_{\overline k}(u)$ and the braid relation for $R$ we have
$$
R_i^{\pm 1}L_{\overline k}(u) = L_{\overline k}(u)R_i^{\pm 1}\qquad \forall\,i\not\in\{k-1,k\}.
$$
Taking into account this property we can write (we  omit the arguments in the matrices)
$$
L_{\overline{k+1}}L_{\overline{k+2}} = R_kL_{\overline{k}}R_k^{-1}R_{k+1}R_kL_{\overline{k}}R_k^{-1}R_{k+1}^{-1} =
R_kR_{k+1}(L_{\overline{k}}L_{\overline{k+1}})R_{k+1}^{-1}R_k^{-1},
$$
where we used the braid relation $R_k^{-1}R_{k+1}R_k = R_{k+1}R_k R_{k+1}^{-1}$ in the middle of the above expression.

Then we continue as follows
$$
L_{\overline{k+1}}L_{\overline{k+2}}L_{\overline{k+3}} = R_kR_{k+1}L_{\overline{k}}L_{\overline{k+1}}L_{\overline{k+3} }R_{k+1}^{-1}R_k^{-1} =
R_kR_{k+1}R_{k+2}(L_{\overline{k}}L_{\overline{k+1}}L_{\overline{k+2}})R_{k+2}^{-1}R_{k+1}^{-1}R_{k}^{-1}
$$
and so on till the following intermediate result:
$$
L_{\overline{k+1}}\dots L_{\overline{k+p}} = R_k\dots R_{k+p-1}(L_{\overline{k}}\dots L_{\overline{k+p-1}}) R_{k+p-1}^{-1}\dots R_{k}^{-1}.
$$

Now, with the use of definition (\ref{a-sym}) and the braid relation for $R$ we can show that
$$
A^{(p)}_{k+1\dots k+p}R_kR_{k+1}\dots R_{k+p-1} = R_kR_{k+1}\dots R_{k+p-1}A^{(p)}_{k\dots k+p-1}.
$$
So, we have
\begin{eqnarray*}
{\rm Tr}_{R(k+1\dots k+p)}\!\!\!\!\!&&\!\!\!\!\! A^{(p)}_{k+1\dots k+p}L_{\overline{k+1}}L_{\overline{k+2}}\dots L_{\overline{k+p}} = \\
&&{\rm Tr}_{R(k+1\dots k+p)}R_kR_{k+1}\dots R_{k+p-1}X_{1\dots k+p-1}R_{k+p-1}^{-1}\dots R_k^{-1},
\end{eqnarray*}
where for the sake of simplicity we denoted $X_{1\dots k+p-1} = A^{(p)}_{k\dots k+p-1}L_{\overline k}\dots L_{\overline{k+p-1}}$. 

At last, applying the relation
(\ref{tr-property}) $p$ times we get
$$
{\rm Tr}_{R(k+1\dots k+p)} A^{(p)}_{k+1\dots k+p}L_{\overline{k+1}}L_{\overline{k+2}}\dots L_{\overline{k+p}} =
I_{1\dots k}{\rm Tr}_{R(k\dots k+p-1)} A^{(p)}_{k\dots k+p-1}L_{\overline{k}}L_{\overline{k+1}}\dots L_{\overline{k+p-1}}.
$$
Thus, we have shifted down by one the numbers of spaces where the $R$-trace is taken. The rest of the proof is evident: applying the same steps
as above, we finally come to the relation (\ref{elem-sym-sh}).\hfill\rule{6.5pt}{6.5pt}

\begin{remark} \rm
Note that the relation  (\ref{elem-sym-sh}) is also valid for a generalized Yangian associated with a compatible pair $(R,F)$, where $R$
is assumed to be a Hecke symmetry. The proof is similar to that above except for the definition of trace: the $R$-trace should be changed for the $F$-trace.
This means that the matrix $C$ in definition (\ref{R-tr}) corresponds to the braiding $F$ instead of $R$.
\end{remark}

\medskip

It is natural to ask whether there exist analogs of other symmetric polynomials  in the algebra $\Y(R)$? We are able to give the positive answer for the so-called
power sums and full symmetric polynomials.  As for the quantum analogs of the Schur polynomials, we do not know their
consistent definition.

\begin{definition} \rm
The elements
\be
p_k(u) = {\rm Tr}_{R(1\dots k)}\Big(L_{\overline 1}(q^{-2(k-1)}u)L_{\overline 2}(q^{-2(k-2)}u)\dots L_{\overline k}(u)R_{k-1}\dots R_2R_1\Big),
\quad k\ge 1.
\label{power-p}
 \ee
are called the {\em quantum power sums}.
\end{definition}

The sets of elementary symmetric polynomials and power sums are related by a series of quantum analogs of the Newton identities \cite{GS1}:
\be
k_qe_k(u) - q^{k-1}p_1(q^{-2(k-1)}u)e_{k-1}(u)+q^{k-2}p_2(q^{-2(k-2)}u)e_{k-2}(u)+\dots +(-1)^kp_k(u)= 0.
\label{q-Newt}
\ee

Similarly to the classical case, these formulae enable us to express power sums via the elementary symmetric polynomials and visa versa. Thus, we can conclude that
the power sums also generate the Bethe subalgebra of $\Y(R)$.

In the rational case we respectively define elementary symmetric polynomials and power sums as follows
\be
e_k(u) = {\rm Tr}_{R(1\dots k)}\Big({A}^{(k)}_{12\dots k}L_{\overline 1}(u)L_{\overline 2}(u-1)\dots L_{\overline k}(u-k+1)\Big),
\label{elem}
\ee
$$
p_k(u) = {\rm Tr}_{R(1\dots k)}\Big(L_{\overline 1}(u-k+1)L_{\overline 2}(u-k+2)\dots L_{\overline k}(u)R_{k-1}R_{k-2}\dots R_1 \Big).
$$

They are related by the following analogs of the Newton identities
$$
k e_k(u) - p_1(u-k+1)e_{k-1}(u)+p_2(u-k+2)e_{k-2}(u)+\dots +(-1)^kp_k(u)= 0.
$$
Also, in the algebras $\Y(R)$ some analogs of the Cayley-Hamilton identity are valid (see \cite{GS1}) but we do not need them here.

\section{Commutativity of Bethe subalgebra}

In this section we prove that the family  of elements $e_k(u)$, $1\le k\le m$, generate a commutative subalgebra in $\Y(R)$. If we expand
each $e_k(u)$ in a series in the inverse powers of $u$, we get countable set of polynomials in entries of the matrices $L[k]$ commuting
with each other.

We begin with some  technical results. Let us introduce a notation for specific products (or chains) of current $R$-matrices:
\be
\Big[R_{i \rightarrow j} (u)\Big]^{(\pm)}  =
 \begin{cases}
R_i (u) R_{i+1} (q^{\pm 2}u) \ldots R_j (q^{\pm 2(j-i)}u), \quad {\rm if} \quad j \ge i  \\
R_i (u) R_{i-1} (q^{\pm 2}u) \ldots R_j (q^{\pm 2(i-j)}u), \quad {\rm if} \quad i\ge j.
 \end{cases}  \label{co1}
\ee
Analogous notation will be applied to the products of $R^{-1}(u)$. Note, that the above chains of current
$R$-matrices are non-singular except for a finite number of isolated values of the parameter $u$ and the inverse chain
is given by the formula:
\begin{equation}
\left\{\Big[R_{i\rightarrow j}(u)\Big]^{(\pm)}\right\}^{-1} = \Big[R^{-1}_{j\rightarrow i}(q^{\pm 2|i-j|}u)\Big]^{(\mp)}.
\label{inv-chain}
\end{equation}

Also, note that the  middle term in formula (\ref{a-sym}) can be presented as follows
$$
\frac{q^k}{k_q}\,I - R_k = - R_k(q^{2k}).
$$
This allows us to get an explicit presentation for the skew-symmetrizer $A^{(k)}$ as a product of chains of $R$-matrices. Namely, we have

\begin{lemma}
\label{lem:2}
The skew-symmetrizers $A^{(k)}$ defined by {\rm (\ref{a-sym})} can be written in the following equivalent forms:\rm
\begin{eqnarray}
A^{(k)}_{1\dots k} \!\!\!&=&\!\!\! \frac{(-1)^{\frac{k(k-1)}{2}}}{k_q!}\,\Big[R_{1\rightarrow k-1}(q^2)\Big]^{(+)}\Big[R_{1\rightarrow k-2}(q^2)\Big]^{(+)}\dots
R_1(q^2),
\label{a-str-1}\\
A^{(k)}_{1\dots k} \!\!\!&=&\!\!\! \frac{(-1)^{\frac{k(k-1)}{2}}}{k_q!}\,\Big[R_{k-1\rightarrow 1}(q^2)\Big]^{(+)}\Big[R_{k-1\rightarrow 2}(q^2)\Big]^{(+)}\dots
R_{k-1}(q^2),
\label{a-str-3}\\
A^{(k)}_{1\dots k} \!\!\!&=&\!\!\!\frac{(-1)^{\frac{k(k-1)}{2}}}{k_q!}\,R_1(q^2)\,\Big[R_{2\rightarrow 1}(q^4)\Big]^{(-)}\dots
\Big[R_{k-1\rightarrow 1}(q^{2(k-1)})\Big]^{(-)},
\label{a-str-2}\\
A^{(k)}_{1\dots k} \!\!\!&=&\!\!\!\frac{(-1)^{\frac{k(k-1)}{2}}}{k_q!}\,R_{k-1}(q^2)\,\Big[R_{k-2\rightarrow k-1}(q^4)\Big]^{(-)}\dots
\Big[R_{k-1\rightarrow 1}(q^{2(k-1)})\Big]^{(-)}.
\label{a-str-4}
\end{eqnarray}
\em
Here, we use the standard notation $k_q! = 1_q2_q\dots k_q$.
\end{lemma}

\medskip

The lemma can be easily proved by induction in $k$ on the base of (\ref{a-sym}).

Then, with the use of the Yang-Baxter equation (\ref{YB}) we can find the following rule for the permutation of $R$-chains:
$$
\Big[R_{1\rightarrow k}(u)\Big]^{(+)}\Big[R_{1\rightarrow k-1}(q^2)\Big]^{(+)} = \Big[R_{2\rightarrow k}(q^2)\Big]^{(+)}
\Big[R_{1\rightarrow k-1}(q^2 u)\Big]^{(+)}R_k(u).
$$
Note, that in the right hand side we have a  cyclic permutation of the set of $u$-depending parameters of the corresponding $R$-chain:
$$
\{u, q^2u,\dots,q^{2(k-1)}u\}\rightarrow \{q^2u,\dots,q^{2(k-1)}u, u\}.
$$
This observation with formulae of Lemma \ref{lem:2} allows us to prove the following result.

\begin{lemma}
\label{lem:3}
The $R$-chains commute with the skew-symmetrizers in accordance with the rule:
\begin{eqnarray}
&&\Big[R_{1\rightarrow k}(q^{-2(k-1)}u)\Big]^{(+)}A_{1\dots k}^{(k)} = A_{2\dots k+1}^{(k)}\Big[R_{1\rightarrow k}(u)\Big]^{(-)},\nonumber\\
\rule{0pt}{5mm}
&&\Big[R_{1\rightarrow k}^{-1}(u)\Big]^{(-)}A_{1\dots k}^{(k)} = A_{2\dots k+1}^{(k)}\Big[R_{1\rightarrow k}^{-1}(q^{-2(k-1)}u)\Big]^{(+)},\nonumber\\
\rule{0pt}{5mm}
&&A^{(k)}_{1\dots k}\Big[R_{k\rightarrow 1}(u)\Big]^{(-)}= \Big[R_{k\rightarrow 1}(q^{-2(k-1)}u)\Big]^{(+)}A^{(k)}_{2\dots k+1},\nonumber\\
&&A^{(k)}_{1\dots k}\Big[R_{k\rightarrow 1}^{-1}(q^{-2(k-1)}u)\Big]^{(+)}= \Big[R_{k\rightarrow 1}^{-1}(u)\Big]^{(-)}A^{(k)}_{2\dots k+1}.\nonumber
\end{eqnarray}
\end{lemma}
Let us point out the inverse order of the parameters of the $R$-chains in the left and right hand sides
of the above relations, as well as the shift of the spaces in which the skew-symmetrizer acts.

The last auxiliary result is connected with the chains of generating matrices $L(u)$. First of all, applying (\ref{sys1-sh}) we get
\begin{eqnarray*}
\Big[R_{k-1\rightarrow1}(q^{2(k-1)})\Big]^{(-)}L_{\overline 1}(u)L_{\overline 2}(q^{-2}u)\!\!\!&\dots&\!\!\!  L_{\overline {k-1}}(q^{-2(k-2)}u)L_{\overline k}(q^{-2(k-1)}u) = \\
L_{\overline 1}(q^{-2}u)L_{\overline 2}(q^{-4}u)\!\!\!&\dots&\!\!\! L_{\overline {k-1}}(q^{-2(k-1)}u)L_{\overline k}(u) \Big[R_{k-1\rightarrow1}(q^{2(k-1)})\Big]^{(-)}.
\end{eqnarray*}
Then, using Lemma \ref{lem:2} we come to the following result:
\be
A^{(k)}_{1\dots k}L_{\overline 1}(u)L_{\overline 2}(q^{-2}u)\dots L_{\overline k}(q^{-2(k-1)}u) =
L_{\overline 1}(q^{-2(k-1)}u)L_{\overline 2}(q^{-2(k-2)}u)\dots L_{\overline k}(u) A^{(k)}_{1\dots k}.
\label{A-L-ch}
\ee

Now, we are ready to prove the main theorem of this paper.

\begin{proposition}
The elementary symmetric functions {\rm (\ref{elem-sym})} commute with each other:
\rm
\be
e_k (u) e_p (v)=e_p (v) e_k (u), \qquad \forall\,k,p\ge 1,\quad\forall\,u,v.
\label{e-comm}
\ee
\end{proposition}
\noindent
{\bf Proof.} In order to simplify formulae  we introduce a shorthand notation for a chain of the  matrices $L(u)$ with shifted arguments:
$$
\big[L_{\overline i\rightarrow \overline j}(u)\big]^{(\pm)} = L_{\overline i}(u)L_{\overline {i+1}}(q^{\pm 2}u)\dots L_{\overline j}(q^{\pm 2(j-i)}u).
$$
Using this notation and the identity (\ref{elem-sym-sh}) we can write the left hand side of (\ref{e-comm}) in the form:
\be
e_k(u)e_p(v) = {\rm Tr}_{R(1\dots k+p)}A^{(k)}_{1\dots k}A^{(p)}_{k+1\dots k+p}\big[L_{\overline{1}\rightarrow \overline{k}}(u)\big]^{(-)}
\big[L_{\overline{k+1}\rightarrow \overline{k+p}}(v)\big]^{(-)}.
\label{start}
\ee

Then, we permute the $u$-depending chain of $L$ matrices with $v$-depending chain with the help of braided Yangian relations (\ref{sys1-sh}).
At the first step of this process we move the utmost right element $L_{\overline k}(q^{-2(k-1)}u)$ through the $v$-chain of the matrices $L$.
Relations (\ref{sys1-sh}) give:
$$
L_{\overline k}(q^{-2(k-1)}u)L_{\overline{k+1}}(v) = R^{-1}_k(q^{-2(k-1)}u/v)L_{\overline k}(v)L_{\overline{k+1}}(q^{-2(k-1)}u)R_{k}(q^{-2(k-1)}u/v).
$$
Repeating this procedure we find:
\begin{eqnarray*}
L_{\overline k}\!\!\!\!\!&&\!\!\!\!\!(q^{-2(k-1)}u)\big[L_{\overline{k+1}\rightarrow \overline{k+p}}(v)\big]^{(-)} =\\
\!\!\!&&\!\!\!\Big[R^{-1}_{k\rightarrow k+p-1}(q^{-2(k-1)}x)\Big]^{(+)}\big[L_{\overline{k}\rightarrow \overline{k+p-1}}(v)\big]^{(-)}L_{\overline {k+p}}(q^{-2(k-1)}u)
\Big[R_{k+p-1\rightarrow k}(q^{-2(k-p)}x)\Big]^{(-)},
\end{eqnarray*}
where $x=u/v$.

Now,  we should move the next factor $L_{\overline{k-1}}(q^{-2(k-2)}u)$ to the right position (it is possible since this factor commutes with the chain of $R^{-1}$-matrices
appeared in the above formula) and so on. Thus, we arrive to the expression:
\begin{eqnarray}
\big[L_{\overline 1\rightarrow \overline k}(u)\big]^{(-)}\big[L_{\overline {k+1}\rightarrow \overline {k+p}}(v)\big]^{(-)}=&&\nonumber \\
\prod_{1\le s\le k}^{\leftarrow}
\Big[R^{-1}_{s\rightarrow s+p-1}\!\!\!\!\!&&\!\!\!\!\!(q^{-2(s-1)}x)\Big]^{(+)}\big[L_{\overline 1\rightarrow \overline p}(v)\big]^{(-)}\big[
L_{\overline {p+1}\rightarrow \overline {p+k}}(u)\big]^{(-)}\nonumber\\
&&\qquad\prod_{1\le r\le k}^{\rightarrow}\Big[R_{r+p-1\rightarrow r}(q^{-2(r-p)}x)\Big]^{(-)},\label{perm-1}
\end{eqnarray}
where we use the following notation for an ordered product of noncommutative factors:
$$
\prod_{1\le i\le k}^{\rightarrow}Q_i \equiv Q_1Q_2\dots Q_k, \qquad  \prod_{1\le i\le k}^{\leftarrow}Q_i \equiv Q_kQ_{k-1}\dots Q_1.
$$

So, with the use of (\ref{perm-1}) we rewrite (\ref{start}) in the following form
\begin{eqnarray}
e_k(u)e_p(v) \!\!\!&=& \!\!\! {\rm Tr}_{R(1\dots k+p)}\Big\{A^{(k)}_{1\dots k}A^{(p)}_{k+1\dots k+p} \prod_{1\le s\le k}^{\leftarrow}
\Big[R^{-1}_{s\rightarrow s+p-1}(q^{-2(s-1)}x)\Big]^{(+)}\big[L_{\overline 1\rightarrow \overline p}(v)\big]^{(-)}\nonumber\\
&&\times \big[L_{\overline {p+1}\rightarrow \overline {p+k}}(u)\big]^{(-)}\prod_{1\le r\le k}^{\rightarrow}
\Big[R_{r+p-1\rightarrow r}(q^{-2(r-p)}x)\Big]^{(-)}\Big\}.\label{prom-2}
\end{eqnarray}
Now, we are going to prove that the ordered products of chains of $R-$ and $R^{-1}$-matrices in the above formula cancel each other under the trace operation. Consequently, the right
hand side turns into the product $e_p(v)e_k(u)$.

Note that except for the simplest case $k=p=1$ this cancellation is not straightforward since in (\ref{prom-2}) there is a pair of skew-symmetrizers between two chains of $R$-matrices.
First, we permute the skew-symmetrizers with the $R$-matrices which leads to the changes of their ordering. Then, we have to make some identical transformations in order to
guarantee that the parameters of the current $R$-matrices and and inverse current $R$-matrices fit each other in a proper way (see (\ref{inv-chain})).

So, we begin with moving the skew-symmetrizers through the chains of $R^{-1}$-matrices with the use of the relations of Lemma \ref{lem:3}. For the skew-symmetrizer $A^{(p)}_{k+1\dots k+p}$
this procedure is straightforward (apply the second line of Lemma \ref{lem:3}):
$$
A^{(p)}_{k+1\dots k+p} \prod_{1\le s\le k}^{\leftarrow} \Big[R^{-1}_{s\rightarrow s+p-1}(q^{-2(s-1)}x)\Big]^{(+)} =
\prod_{1\le s\le k}^{\leftarrow} \Big[R^{-1}_{s\rightarrow s+p-1}(q^{-2(s-p)}x)\Big]^{(-)} A^{(p)}_{1\dots p}.
$$
In order to move the second skew-symmetrizer $A^{(k)}_{1\dots k}$ we have to reorder the chains of $R$-matrices using the following evident identity:
$$
\prod_{1\le s\le k}^{\leftarrow} \Big[R^{-1}_{s\rightarrow s+p-1}(q^{-2(s-p)}x)\Big]^{(-)} =
\prod_{1\le s\le p}^{\rightarrow} \Big[R^{-1}_{s+k-1\rightarrow s}(q^{-2(s+k-p-1)}x)\Big]^{(+)}.
$$
Now, we can move $A^{(k)}_{1\dots k}$ with the help of the last line formula of Lemma \ref{lem:3}:
$$
A^{(k)}_{1\dots k}\prod_{1\le s\le p}^{\rightarrow} \Big[R^{-1}_{s+k-1\rightarrow s}(q^{-2(s+k-p-1)}x)\Big]^{(+)} =
\prod_{1\le s\le p}^{\rightarrow} \Big[R^{-1}_{s+k-1\rightarrow s}(q^{-2(s-p)}x)\Big]^{(-)}A^{(k)}_{p+1\dots k+p}.
$$
So, we come to the relation:
\begin{eqnarray*}
e_k(u)e_p(v) \!\!\!&=& \!\!\! {\rm Tr}_{R(1\dots k+p)}\Big\{\prod_{1\le s\le k}^{\rightarrow} \Big[R^{-1}_{k+1-s\rightarrow k+p-s}(q^{-2(s-p)}x)\Big]^{(-)}
A^{(p)}_{1\dots p}A^{(k)}_{p+1\dots k+p} \big[L_{\overline 1\rightarrow \overline p}(v)\big]^{(-)}\\
&&\times \big[L_{\overline {p+1}\rightarrow \overline {p+k}}(u)\big]^{(-)}\prod_{1\le r\le k}^{\rightarrow}
\Big[R_{r+p-1\rightarrow r}(q^{-2(r-p)}x)\Big]^{(-)}\Big\},
\end{eqnarray*}
where  we also reordered the chains of inverse $R$-matrices with the use of identity:
$$
\prod_{1\le s\le p}^{\rightarrow} \Big[R^{-1}_{s+k-1\rightarrow s}(q^{-2(s-p)}x)\Big]^{(-)} =
\prod_{1\le s\le k}^{\rightarrow} \Big[R^{-1}_{k+1-s\rightarrow k+p-s}(q^{-2(s-p)}x)\Big]^{(-)}.
$$
Next, using the projector property (\ref{a-sym-pro}) we can write
$$
A_{1\dots p}^{(p)}A_{p+1\dots k+p}^{(k)} = A_{1\dots p}^{(p)}A_{p+1\dots k+p}^{(k)}A_{1\dots p}^{(p)}A_{p+1\dots k+p}^{(k)}
$$
and then  move the second pair of skew-symmetrizers to the rightmost position. In so doing, we apply (\ref{A-L-ch}) to pass through the chains of $L$-matrices and
Lemma \ref{lem:3} to permute $A^{(k)}A^{(p)}$ with the chains of current $R$-matrices.

At last, with the help of (\ref{cyclic-p}) we cyclically move the chains of $R$-matrices and the second pair of skew-symmetrizers. Finally, we get the following result:
\begin{eqnarray*}
e_k(u)e_p(v) \!\!\!&=& \!\!\! {\rm Tr}_{R(1\dots k+p)}\Big\{\prod_{1\le r\le k}^{\rightarrow} \Big[R_{r+p-1\rightarrow r}(q^{-2(k-r)}x)\Big]^{(+)}
A^{(k)}_{1\dots k}A^{(p)}_{k+1\dots k+p}\\
\!\!\!&\times&\!\!\!\!\!\prod_{1\le s\le k}^{\leftarrow} \Big[R^{-1}_{s\rightarrow s+p-1}(q^{-2(k-p+1-s)}x)\Big]^{(-)}A^{(k)}_{p+1\dots k+p}A^{(p)}_{1\dots p}\\
\!\!\!&\times&\!\!\!\!\! \big[L_{\overline 1\rightarrow \overline p}(q^{-2(p-1)}v)\big]^{(+)}\big[L_{\overline {p+1}\rightarrow \overline {p+k}}(q^{-2(k-1)}u)\big]^{(+)}\Big\}.
\end{eqnarray*}

Note, that in the above expression the chains of current $R$-matrices are exactly inverse to each other (see (\ref{inv-chain})). The only obstacle for their cancellation is the
pair of skew-symmetrizers  between them. To remove these skew-symmetrizers,  we permute the {\it outer} pair $A^{(k)}_{p+1\dots k+p}A^{(p)}_{1\dots p}$ with the chain
of $R^{-1}$-matrices, use the above projector property and then move the pair of skew-symmetrizers back through the chain of inverse $R$-matrices.

Due to such a twice permutation the chains of inverse current $R$-matrices are not changed and all $R$-matrix chains  cancel each other. Finally, we get
\begin{eqnarray*}
e_k(u)e_p(v) \!\!\!&=& \!\!\! {\rm Tr}_{R(1\dots k+p)}\Big\{A^{(p)}_{1\dots p}A^{(k)}_{p+1\dots k+p}\big[L_{\overline 1\rightarrow \overline p}(q^{-2(p-1)}v)\big]^{(+)}
\big[L_{\overline {p+1}\rightarrow \overline {p+k}}(q^{-2(k-1)}u)\big]^{(+)}\Big\} \\
\!\!\!&=& \!\!\! {\rm Tr}_{R(1\dots k+p)}\Big\{\big[L_{\overline 1\rightarrow \overline p}(v)\big]^{(-)}
\big[L_{\overline {p+1}\rightarrow \overline {p+k}}(u)\big]^{(-)}A^{(p)}_{1\dots p}A^{(k)}_{p+1\dots k+p}
\Big\}= e_p(v)e_k(u).
\end{eqnarray*}

As was noticed above, on expanding $e_k(u)$ and $e_p(v)$ in a series in $u^{-1}$ and $v^{-1}$ we can reformulate this commutativity in terms of plynomials in
entries of the matrices $L[k]$.  \hfill \rule{6.5pt}{6.5pt}

\section{Braided Gaudin models}

In the present section we deal with some versions of the Gaudin model, associated with the braided Yangians considered above. First, we recall the classical version of the rational
Gaudin model.

Let $M(k)$, $k=1,2,\dots,K$, be $m\times m$ matrices subject to the following relations
\be
M(k)_1 M(l)_2= M(l)_2  M(k)_1,\quad  k\not=l,
\label{first}
\ee
and
\be
M(k)_1 M(k)_2- M(k)_2 M(k)_1=M(k)_1 P-M(k)_2 P.
\label{sec}
\ee
Then the elements
\be
H_k= \sum_{l\not=k} {\rm Tr} \,\frac{M(k) M(l)}{u_k-u_l},
\label{ham}
\ee
where $u_k,$ $k=1,2,\dots,K$, are pairwise distinct complex numbers, commute with each other.

This claim results from the fact that the elements $H_k$ are  quadratic Hamiltonians of the rational Gaudin model. The higher Hamiltonians of this model have been
constructed by D.Talalaev \cite{T}. His result can be presented as follows.

Let $L(u)=\sum_{k\ge 1}L[k]u^{-k}$ be an $m\times m$ matrix subject to the relation
\be
[L_1(u), L_2(v)]=\Big[\frac{P}{u-v}\,, L_1(u)+ L_2(v)\Big].
\label{PP}
\ee
Then the elements
\be
QH_k(u)=\Tr_{(1\dots m)}A^{(m)}_{1\dots m} (L_{ 1}(u)-I\, \frac{d}{d\, u})(L_{ 2}(u)-I\, \frac{d}{d\, u})\dots (L_{ k}(u)-I\, \frac{d}{d\, u})\triangleright 1 \label{QH}
\ee
commute with each other:
$$
QH_k(u) QH_l(v)=QH_l(v) QH_k(u)\quad \forall \, u, v \quad {\rm and}\quad  \forall\, k, l.
$$
Here the notation $\triangleright$ means that  the above differential operator is applied  to the unit. 
Thus, all terms of the resulting differential operator $\sum_{p}^k\, Q_p(u) \pa^p$ are
cancelled except that which does not contain the derivative.

Now, by using the fact that the matrix $\sum_{i=1}^{K} \frac{M_i}{u-u_i}$ satisfies the relation (\ref{PP}) (a proof of this fact is straightforward)
we get a family of Hamiltonians for the rational Gaudin model.

Observe that the system (\ref{sec})  means that entries of each matrix $M(k)$ generate a Lie algebra $gl(m)$, while the entries of matrices $M(k)$ and $M(l)$ with $k\not= l$
commute with each other due to (\ref{first}). Thus, the vector space spanned by
the entries of all matrices $M(k)$ is endowed with the structure of the Lie algebra $gl(m)^{\oplus K}$.

Now, in the relations (\ref{first}) and (\ref{sec}) we replace all indices by their overlined counterparts (see definition (\ref{transf})) assuming $R$ to be an involutive symmetry.
Namely, we consider the following relations
\be
M(k)_{\ov 1} M(l)_{\ov 2}= M(l)_{\ov 2} M(k)_{\ov 1},\quad  k\not=l,
\label{first1}
\ee
and
\be
M(k)_{\ov 1} M(k)_{\ov 2}- M(k)_{\ov 2} M(k)_{\ov 1}=M(k)_{\ov 1} R-M(k)_{\ov 2} R.
\label{sec1}
\ee

Note that the relations (\ref{sec1}) define the enveloping algebra of a braided Lie algebra generated by entries of the matrix $M(k)$ as defined in \cite{GPS}.
Our present purpose is to show that the elements
\be
H_k= \sum_{l\not=k} {\rm Tr}_R \,\frac{M(k)\, M(l)}{u_k-u_l}
\label{ham1}
\ee
commute with each other. Emphasize, that in contrast with (\ref{ham}) the trace in this formula is braided.

In order to prove the above claim we follow the main lines of the Talalaev's approach but instead of the Drinfeld's Yangian we consider the braided Yangian arising from an
involutive symmetry $R$ of bi-rank $(m|0)$.

More precisely, we work with the braided Yangian defined by (\ref{sys1}) with the following current $R$ matrix:
$$
R(u,v)=R-\frac{h\, I}{u-v},
$$
which differs from the initial rational current $R$-matrix by the numerical parameter $h$ in the numerator. We denote the corresponding braided Yangian $\Y(R)_h$.

Also, we pass to a new generating matrix $\tilde L(u)= \|\tilde{l}_i^j(u)\|_{1\leq i, j \leq m}$ connected with the initial one by a linear shift:
$$
L(u)=I+h\,\tilde{L}(u) :\quad l_i^j(u)=\de_i^j+h\,\tilde{l}_i^j(u).
$$
Consequently, we get the following system for the generating matrix $\tilde{L}(u)$:
 $$
(R-\frac{h\, I}{u-v}) {\tL(u)}_1 R {\tL(v)}_1-{\tL(v)}_1 R {\tL(v)}_1 (R-\frac{h\, I}{u-v} )=
 -[R, \frac{\tL_1(u)-\tL_1(v)}{u-v}].
$$

In the limit $h\to 0$ we get an algebra defined by the following system
\be [\tL_{\ov 1}(u), \tL_{\ov 2}(v)]=\Big[\frac{R}{u-v}, \tL_{\ov 1}(u)+\tL_{\ov 2}(v)\Big].
\label{skl}
\ee

This relation is an analog of the first Sklyanin structure and  can be treated in terms of the braided Lie algebras (see \cite{GPS}).

In the algebra $\Y(R)_h$ we consider the following elements
\be
\hat{e}_k (u) = {\rm Tr}_{R(12 \ldots m)} \,A^{(m)}_{12 \dots m} {L}_{\overline{1}} (u) {L}_{\overline{2}} (u-h) \ldots {L}_{\overline{k}} (u-h(k-1)), \quad k\ge 1.
\label{he}
\ee
For $h=1$ the element $ \hat{e}_k (u)$ differs from ${e}_k (u)$ defined in (\ref{elem}) by a non-trivial multiplier\footnote{This property is only valid in the braided
 Yangians of RE type.} (see \cite{GS3}) and consequently we can conclude that the elements $\hat{e}_k (u)$  commute with each other in the algebra $\Y(R)_h$.

Next, we expand the elements $\hat{e}_k (u)$ in a series in $h$ and introduce the following linear combinations
$$
\tau_k(u)=\sum_{p=0}^k(-1)^{k-p}\frac{k!}{p! (k-p)!}\,\hat{e}_p(u).
$$
The corresponding $h$-expansion of $\tau_k(u)$ begins with a term proportional to $h^k$ (see \cite{T,  GS2}). So, we arrive to the following conclusion:
the elements
\be
QH_k(u)=\Tr_{R(1\dots m)}A^{(m)}_{1\dots m} (\tL_{\ov 1}(u)-I\, \frac{d}{d\, u})(\tL_{\ov 2}(u)-I\, \frac{d}{d\, u})\dots (\tL_{\ov k}(u)-I\, \frac{d}{d\, u})\triangleright 1 \label{QH1}
\ee
commute with each other in the algebra defined by (\ref{skl}).

Also, observe that the matrix
$$
 \sum_{k=1}^K \frac{M(k)}{u-u_k},
$$
satisfies the system (\ref{skl}), provided that the matrices $M(k)$ are subject to (\ref{first1}) and (\ref{sec1}). A verification of this property is straightforward.

Consequently, the mutual commutativity of the elements (\ref{QH1}) remains valid if we replace  the matrix $\tL_1(u)$ by $\sum_k \frac{M(k)}{u-u_k}$ in this formula.
As usual, the commutativity of the quadratic Hamiltonians can be deduced from this claim via computing the residues  at the points $u_k$.

Unfortunately, this method does not lead to new interesting systems when it is applied to braided Yangians related to Hecke symmetries. Indeed, in \cite{GS3} a similar
linear shift of the generating matrix $L(u)$ was applied in the braided Yangian $\Y(R)$ corresponding to a Hecke symmetry $R$. However, the role of the parameter $h$
in a subsequent expansion was attributed to  $\log  q$. While $q$ tends to 1 the Hecke symmetry $R$ tends to either the flip $P$ or an involutive symmetry $\tilde R$ if $R$ is
a deformation of $P$ or $\tilde R$ respectively.

Finally,  we come to a limit algebra defined by the system
\be 
[\tL_{\ov 1}(u),\, \tL_{\ov 2}(v)]= \Big[\frac{\tilde R}{u-v}, u\tL_{\ov 1}(u)+v\tL_{\ov 2}(v)\Big], \label{last}
\ee
which does not depend on $q$ (in the former case ${\tilde R}=P$).  

The case ${\tilde R}=P$ was considered in \cite{GS3}. Namely, by using the Talalaev's method, it is easy to
 construct a Bethe subalgebra in  the algebra defined by (\ref{last}). Then, upon applying the morphism
$$
\tL(u)\mapsto  \sum_{k=1}^K \frac{M(k)\, u_k}{u-u_k},
$$
it is possible to conclude that  the elements
\be 
H_k= \sum_{l\not=k} {\rm Tr}\,\frac{M(k)\, M(l)\, u_l}{u_k-u_l} \label{lastt}
\ee
commute with each other, provided the matrices $M(k)$ are subject to (\ref{first}) and (\ref{sec}).
However, this model is not new, since the change of parameters $u_k\mapsto u_k^{-1}$ reduces this family to that (\ref{ham}) above. 

In general, i.e. while ${\tilde R}\not=P$, all these considerations remain valid if we assume that
the matrices $M(k)$ are subject to (\ref{first1}) and (\ref{sec1}) with $R={\tilde R}$    and replace
 the usual trace $\rm{\Tr}$ in (\ref{lastt}) by its ${\tilde R}$-counterpart.

 \end{document}